\documentclass{amsart}

\usepackage{amssymb, amsmath}
\usepackage{mathrsfs}
\usepackage{amscd}
\usepackage{verbatim}

\usepackage[colorlinks,linkcolor={blue},
citecolor={blue},urlcolor={red},]{hyperref}

\theoremstyle{plain}
\newtheorem{theorem}{Theorem}[section]
\theoremstyle{remark}
\newtheorem{remark}{Remark}[section]
\newtheorem{example}{Example}[section]
\theoremstyle{plain}
\newtheorem{corollary}{Corollary}[section]
\newtheorem{lemma}{Lemma}[section]
\newtheorem{proposition}{Proposition}[section]

\numberwithin{equation}{section}

\def\Z{{\mathbb Z}}

\def\R{{\mathbb R}}

\newcommand{\E}{{\mathbb E}}
\renewcommand{\P}{{\mathbb P}}
\newcommand{\F}{{\mathcal F}}

\newcommand{\g}{\gamma}

\newcommand{\e}{\varepsilon}

\renewcommand{\O}{\Omega}

\newcommand{\beq}{\begin{equation}}
\newcommand{\eeq}{\end{equation}}
\newcommand{\bal}{\begin{aligned}}
\newcommand{\eal}{\end{aligned}}
\newcommand{\ben}{\begin{enumerate}}
\newcommand{\een}{\end{enumerate}}
\newcommand{\bit}{\begin{itemize}}
\newcommand{\eit}{\end{itemize}}

\newcommand{\bth}{\begin{theorem}}
\renewcommand{\eth}{\end{theorem}}
\newcommand{\bpr}{\begin{proposition}}
\newcommand{\epr}{\end{proposition}}
\newcommand{\ble}{\begin{lemma}}
\newcommand{\ele}{\end{lemma}}
\newcommand{\bpf}{\begin{proof}}
\newcommand{\epf}{\end{proof}}
\newcommand{\bex}{\begin{example}}
\newcommand{\eex}{\end{example}}
\newcommand{\bre}{\begin{example}}
\newcommand{\ere}{\end{example}}

\newcommand{\calL}{{\mathcal L}}

\newcommand{\one}{{{\bf 1}}}
\newcommand{\embed}{\hookrightarrow}
\newcommand{\s}{^*}
\newcommand{\lb}{\langle}
\newcommand{\rb}{\rangle}

\begin{document}

\title
[Regularity of Brownian motions in infinite dimensions]{On Besov
regularity of Brownian motions in infinite dimensions}

\author{Tuomas Hyt\"onen}
\address{Department of Mathematics and Statistics\\
University of Helsinki\\ Gustaf H\"allstr\"omin katu 2B \\ FI-00014
Helsinki\\ Finland} \email{tuomas.hytonen@helsinki.fi}
\thanks{T.~Hyt\"onen is supported by the Academy of Finland (grant 114374).}

\author{Mark Veraar}

\address{Institute of Mathematics\\
Polish Academy of Sciences \\ Sniadeckich 8 \\ 00-950 Warsaw
\\ Poland} \email{m.veraar@impan.gov.pl, mark@profsonline.nl}
\thanks{M.~C.~Veraar is supported by the Netherlands Organisation for Scientific
Research (NWO) 639.032.201 and by the Research Training Network
MRTN-CT-2004-511953}


%
%
%
%

\maketitle
\begin{abstract}
We extend to the vector-valued situation some earlier work of
Ciesielski and Roynette on the Besov regularity of the paths of the
classical Brownian motion. We also consider a Brownian motion as a
Besov space valued random variable. It turns out that a Brownian
motion, in this interpretation, is a Gaussian random variable with
some pathological properties. We prove estimates for the first
moment of the Besov norm of a Brownian motion. To obtain such
results we estimate expressions of the form $\E \sup_{n\geq
1}\|\xi_n\|$, where the $\xi_n$ are independent centered Gaussian
random variables with values in a Banach space. Using isoperimetric
inequalities we obtain two-sided inequalities in terms of the first
moments and the weak variances of $\xi_n$.
\end{abstract}

\section{Introduction}\label{sec:intro}

Let $(\O, \mathcal A, \P)$ be a complete probability space. Let
$W:[0,1]\times\O\to \R$ be a standard Brownian motion. Since $W$ has
continuous paths, it is easy to check that $W:\O\to C([0,1])$ is a
$C([0,1])$-valued Gaussian random variable. Moreover, since $W$ is
$\alpha$-H\"older continuous for all $\alpha\in (0,\frac12)$ one can
also show that for all $0<\alpha<1/2$, $W:\O\to C^{\alpha}([0,1])$
is a Gaussian random variable. In this way one obtains results like
\[\E \exp(\e\|W\|_{C^{\alpha}([0,1])}^2) < \infty\]
for some $\e>0$.

In \cite{Cie,Cie2} Ciesielski has improved the H\"older continuity
results of Brownian motion using Besov spaces. He has proved that
almost all paths of $W$ are in the Besov space
$B^{1/2}_{p,\infty}(0,1)$ for all $p\in [1, \infty)$ or even in the
Besov-Orlicz space $B^{1/2}_{\Phi_2,\infty}(0,1)$, where $\Phi_2(x)
= e^{x^2}-1$ (for the definition we refer to Section
\ref{sec:prel}). In \cite{Roy} Roynette has characterized the set of
indices $\alpha, p,q$ for which the paths of Brownian motion belong
the Besov spaces $B^{\alpha}_{p,q}(0,1)$.

The proofs of the above results are based on certain coordinate
expansions of
the Brownian motion and 
descriptions of the Besov norms in terms of the corresponding
expansion coefficients of a function. We will give more direct
proofs of these results which employ the usual modulus-of-continuity
definition of the Besov norms. Our methods also carry over to the
vector-valued situation.

\medskip

Let $X$ be a real Banach space. We will write $a\lesssim b$ if there
exists a universal constant $C>0$ such that $a\leq Cb$, and $a\eqsim
b$ if $a\lesssim b\lesssim a$. If the constant $C$ is allowed to
depend on some parameter $t$, we write $a\lesssim_t b$ and
$a\eqsim_t b$ instead. Let $(l^\Theta, \|\cdot\|_{\Theta})$ denote
the Orlicz sequence space with $\Theta(x) = x^2 e^{-\frac{1}{x^2}}$.
Let $(\xi_n)_{n\geq 1}$ be independent centered $X$-valued Gaussian
random variables with weak variances $(\sigma_n)_{n\geq 1}$ and $m =
\sup_{n\geq 1} \E\|\xi_n\|$. In Section \ref{sec:linfty} we will
show that
\begin{equation}\label{eq:estlinfty}
\E\sup_{n\geq 1} \|\xi_n\| \eqsim  m  + \|(\sigma_n)_{n\geq
1}\|_{\Theta}.
\end{equation}
As a consequence of the Kahane--Khinchine inequalities a similar
estimate holds for $(\E\sup_{n\geq 1} \|\xi_n\|^p)^{1/p}$ for all
$p\in [1, \infty)$ as well, at the cost of replacing $\eqsim$ by
$\eqsim_p$.
The proof of \eqref{eq:estlinfty} is based on isoperimetric
inequalities for Gaussian random variables (cf. \cite{LeTa}).

In Section \ref{sec:regBM} we obtain regularity properties of
$X$-valued Brownian motions $W$. In particular we show that for the
paths of an $X$-valued Brownian motion $W$ we have $W\in
B^{1/2}_{p,\infty}(0,1;X)$ for all $p\in [1, \infty)$ or even $W\in
B^{1/2}_{\Phi_2,\infty}(0,1;X)$. Thus we can consider the mappings
$W:\O\to B^{1/2}_{p,\infty}(0,1;X)$ and $W:\O\to
B^{1/2}_{\Phi_2,\infty}(0,1;X)$. A natural question is whether $W$
is a Gaussian random variable with values in one of these spaces. To
answer this some problems have to be solved, because the Banach
spaces $B^{1/2}_{p,\infty}(0,1)$ and $B^{1/2}_{\Phi_2,\infty}(0,1)$
are non-separable. It will be shown in Section~\ref{sec:rvW} that
$W$ is indeed a Gaussian random variable, but it has some peculiar
properties. For instance we find that there exists an $\e>0$ such
that
\[
  \P(\|W\|_{B^{1/2}_{p,\infty}(0,1;X)}\leq \e)=\P(\|W\|_{B^{1/2}_{\Phi_2,\infty}(0,1;X)}\leq \e)=0
\]
which is rather counterintuitive for a centered Gaussian random
variable. It implies in particular that $W$ is not Radon. In the
last Section~\ref{sec:moment} we apply the results from Section
\ref{sec:linfty} to obtain explicit estimates for
$\E\|W\|_{B^{1/2}_{p,\infty}(0,1;X)}$ and
$\E\|W\|_{B^{1/2}_{\Phi_2,\infty}(0,1;X)}$.

%


\section{Preliminaries}\label{sec:prel}

\subsection{Orlicz spaces}\label{subsec:Orlicz}
We briefly recall the definition of Orlicz spaces. More details can
be found in \cite{KrRu,RR,Za}.

Let $(S, \Sigma, \mu)$ be a $\sigma$-finite measure space and let
$X$ be a Banach space. Let $\Phi:\R\to \R_+$ be an even convex
function such that $\Phi(0) = 0$ and $\lim_{x\to \infty}\Phi(x) =
\infty$. The Orlicz space $L^\Phi(S;X)$ is defined as the set of all
strongly measurable functions $f:S\to X$ (identifying functions
which are equal $\mu$-a.e.) with the property that there exists a
$\delta>0$ such that
\[M_{\Phi}(f/\delta) := \int_S\Phi(\|f(s)\|/\delta) \, d\mu(s)<\infty.\]
This space is a vector space and we define
\[\rho_\Phi(f) = \inf\{\delta>0:M_{\Phi}(f/\delta)\leq 1\}.\]
The mapping $\rho_\Phi$ defines a norm on $L^\Phi(S;X)$ and it turns
$L^\Phi(S;X)$ into a Banach space. It is usually referred to as the
\emph{Luxemburg norm}.

For $f\in L^\Phi(S;X)$ we also define the \emph{Orlicz norm}
\[\|f\|_{\Phi} = \inf_{\delta>0}\Big\{\frac{1}{\delta}(1+M_{\Phi}(\delta f))\Big\}.\]
Usually the Orlicz norm is defined in a different way using duality,
but the above norm gives exactly the same number (cf. \cite[Theorem
III.13]{RR}).
%

The two norms are equivalent, as shown in the following:

\begin{lemma}\label{lem:OrliczLux}
For all $f\in L^\Phi(S;X)$ we have
\[\rho_\Phi(f) \leq \|f\|_{\Phi} \leq 2 \rho_\Phi(f).\]
\end{lemma}

\begin{proof}
Let $\delta>0$ be such that $M_{\Phi}(f \delta)\leq 1$. Then
\[\frac{1}{\delta}(1+M_{\Phi}(\delta f)) \leq \frac2\delta.\]
Taking the infimum over all $\delta>0$ such that $M_{\Phi}(f
\delta)\leq 1$ gives the second inequality.

For the first inequality, choose $\alpha>\|f\|_{\Phi}$. Then there
exists a $\delta>0$ such that
\[
  \frac{1}{\delta}(1+M_{\Phi}(\delta f)) \leq \alpha.
\]
Since $\Phi(0) = 0$ and $\Phi$ is convex it holds that
$\Phi(x/\beta)\leq \Phi(x)/\beta$ for all $x\in \R$ and $\beta\geq
1$. Noting that $\alpha\delta\geq 1$ it follows that
\[
  M_\Phi(f/\alpha) = M_\Phi\Big(\frac{\delta f}{\delta\alpha}\Big)
  \leq\frac{M_\Phi(\delta f)}{{\delta\alpha}}\leq 1.
\]
Since $\rho_{\Phi}(f)$ is the infimum over all $\alpha>0$ for which
the previous inequality holds, and it holds for every
$\alpha>\|f\|_{\Phi}$, we conclude that
$\rho_{\Phi}(f)\leq\|f\|_{\Phi}$.
\end{proof}

It is clear from the proof that the lemma holds for all functions
$\Phi:\R_+\to \R$ that satisfy $\Phi(0)=0$ and $\Phi(x/\beta)\leq
\Phi(x)/\beta$ for all $x\in \R_+$ and $\beta\geq 1$. An interesting
example of a non-convex function that satisfies the above properties
is $\Phi(x) = x e^{-1/x^2}$.

\subsection{The Orlicz sequence space $l^{\Theta}$}
We next present a particular Orlicz space which plays an important
role in our studies. The underlying measure space is now $\Z_+$ with
the counting measure, and we consider the function $\Theta:\R\to
\R_+$ defined by
\begin{equation}\label{Phiexp}
\Theta(x) = x^2\exp\big(-\frac{1}{2x^2}\big).
\end{equation}
This function satisfies the assumptions in Subsection
\ref{subsec:Orlicz} and we can associate an Orlicz sequence space
$l^\Theta$ to it. Thus $l^\Theta$ consists of all sequences
$a:=(a_n)_{n\geq 1}$ for which
\[\rho_{\Theta}(a) := \inf\Big\{\delta>0: \sum_{n\geq 1} \frac{a_n^2}{\delta^2}
\exp\Big(-\frac{\delta^2}{2a_n^2}\Big)\leq 1\Big\}<\infty.\]

The following example illustrates the behaviour of
$\rho_{\Theta}(a)$, but also plays a role later on.

\begin{example}\label{ex:geometricSeq}
If $a_n=\alpha^n$, where $\alpha\in[1/2,1)$, then
$\rho_{\Theta}(a)\eqsim\sqrt{\log(1-\alpha)^{-1}}$.
\end{example}

This may be compared with $\|a\|_p\eqsim(1-\alpha)^{-1/p}$, again
for $\alpha\in[1/2,1)$, and $p\in[1,\infty]$.

\begin{proof}
We consider the equivalent Orlicz norm $\|a\|_{\Theta}$. On the one
hand,
\[\begin{split}
  \sum_{n\geq 1}\lambda^2\alpha^{2n}\exp(-\frac{1}{2\lambda^2\alpha^{2n}})
  &\leq\sum_{n\geq 1}\lambda^2\alpha^{2n}\exp(-\frac{1}{2\lambda^2\alpha^2})
  \\ &=\frac{\lambda^2\alpha^2}{1-\alpha^2}\exp(-\frac{1}{2\lambda^2\alpha^2})\\
  &\leq\frac{\lambda^2}{1-\alpha}\exp(-\frac{1}{2\lambda^2}).
\end{split}\]

On the other hand, let $N\in\Z_+$ be such that $\alpha^{2N}\leq
1/2<\alpha^{2(N-1)}$. Then
\[\begin{split}
  \sum_{n\geq 1}\lambda^2\alpha^{2n}\exp(-\frac{1}{2\lambda^2\alpha^{2n}})
  &\geq\sum_{n=1}^N\lambda^2\alpha^{2n}\exp(-\frac{1}{2\lambda^2\alpha^{2N}})\\
  &\geq\lambda^2\alpha^2\frac{1-\alpha^{2N}}{1-\alpha^2}\exp(-\frac{1}{\lambda^2\alpha^2})\\
  &\geq\frac{\lambda^2}{12(1-\alpha)}\exp(-\frac{4}{\lambda^2}).
\end{split}\]

From these observations it follows that
\[
  \|a\|_{\Theta}
  =\inf_{\lambda>0}\frac{1}{\lambda}(1+M_{\Theta}(\lambda a))
  \eqsim\inf_{\lambda>0}\frac{1}{\lambda}(1+\frac{\lambda^2}{1-\alpha}e^{-1/2\lambda^2})
  =:\inf_{\lambda>0}F(\lambda).
\]
The differentiable function $F$ tends to $\infty$ as $\lambda\to 0$
or $\lambda\to\infty$, so its infimum is attained at a point where
$F'(\lambda)=0$. Since
\[
  F'(\lambda)=-\lambda^{-2}+(1-\alpha)^{-1}e^{-1/2\lambda^2}+(1-\alpha)^{-1}e^{-1/2\lambda^2}\lambda^{-2},
\]
where the middle-term is always positive, $F'(\lambda)=0$ can only
happen if
\[
  (1-\alpha)^{-1}e^{-1/2\lambda^2}\leq 1\qquad\text{i.e.},\qquad
  \lambda^{-1}\geq\lambda_0^{-1}:=\sqrt{2\log(1-\alpha)^{-1}}.
\]
But $1/\lambda$ is the first term in $F(\lambda)$, so we have proved
that $F(\lambda)\gtrsim\sqrt{\log(1-\alpha)^{-1}}$ whenever
$0<\lambda\leq\lambda_0$, and moreover there holds
$F(\lambda_0)\eqsim\sqrt{\log(1-\alpha)^{-1}}$, which completes the
proof.
\end{proof}

\subsection{Besov spaces}
We recall the definition of the vector-valued Besov spaces. For the
real case we refer to \cite{Tr1} and for the vector-valued Besov
space we will give the treatise from \cite{Ko}.

Let $X$ be a real Banach space and let $I=(0,1)$. For $\alpha\in
(0,1)$, $p,q\in [1, \infty]$ the {\em vector-valued Besov space}
$B^{\alpha}_{p,q}(I;X)$ is defined as the space of all functions
$f\in L^p(I;X)$ for which the seminorm (with the usual modification
for $q=\infty$)
\[
  \Big(\int\limits_0^1 (t^{-\alpha} \omega_p(f,t))^q \, \frac{dt}{t}\Big)^{1/q}
\]
is finite. Here
\[
  \omega_p(f,t) = \sup_{|h|\leq t} \|s\mapsto f(s+h)-f(s)\|_{L^p(I(h);X)}
\]
with $I(h) = \{s\in I: s+h\in I\}$. The sum of the $L^p$-norm and
this seminorm turn $B^{\alpha}_{p,q}(I;X)$ into a Banach space.
 By a dyadic approximation argument (see \cite[Corollary 3.b.9]{Ko}) one can show
that the above seminorm is equivalent to
\[
  \|f\|_{p,q, \alpha}
  := \Big(\sum_{n\geq 0} \big(2^{n\alpha} \|s\mapsto f(s+2^{-n})-f(s)\|_{L^p(I(2^{-n});X)}\big)^{q}\Big)^{1/q}
\]
For the purposes below it will be convenient to take
\[
  \|f\|_{B^{\alpha}_{p,q}(I;X)} = \|f\|_{L^p(I;X)}+\|f\|_{p,q, \alpha}
\]
as a Banach space norm on $B^{\alpha}_{p,q}(I;X)$.

For $0<\beta<\infty$, we also introduce the exponential Orlicz and
Orlicz--Besov (semi)norms
\begin{equation*}\begin{split}
  \|f\|_{\mathfrak{L}^{\Phi_{\beta}}(I;X)}
  &:=\sup_{p\geq 1} p^{-1/\beta}\|f\|_{L^p(I;X)},\\
  \|f\|_{\Phi_{\beta},\infty,\alpha}
  &:=\sup_{n\geq 1}2^{\alpha n}\|f-f(\cdot-2^{-n})\|_{\mathfrak{L}^{\Phi_{\beta}}(I(2^{-n});X)}
  =\sup_{p\geq 1}p^{-1/\beta}\|f\|_{p,\infty,\alpha},
\end{split}\end{equation*}
and finally the Orlicz--Besov norm
\begin{equation*}
  \|f\|_{B^{\alpha}_{\Phi_{\beta},\infty}(I;X)}
  :=\sup_{p\geq 1}p^{-1/\beta}\|f\|_{B^{\alpha}_{p,\infty}(I;X)}
  \eqsim \|f\|_{\mathfrak{L}^{\Phi_{\beta}}(I;X)}
  +\|f\|_{\Phi_{\beta},\infty,\alpha}.
\end{equation*}
Because of the inequalities between different $L^p$ norms, it is
immediate that we have equivalent norms above, whether we understand
$p\geq 1$ as $p\in[1,\infty)$ or $p\in\{1,2,\ldots\}$. For
definiteness and later convenience, we choose the latter.

The above-given norm of $\mathfrak{L}^{\Phi_{\beta}}(I;X)$ is
equivalent to the usual norm of the Orlicz space
$L^{\Phi_{\beta}}(I;X)$ from Subsection \ref{subsec:Orlicz} where
$\Phi_{\beta}(x)=\exp(|x|^{\beta})-1$ for $\beta\geq 1$. For
$0<\beta<1$, the function $\Phi_{\beta}$ must be defined in a
slightly different way, but it is still essentially
$\exp(|x|^{\beta})$; see \cite{Cie2}.

For $\beta\in \Z_+\setminus\{0\}$ one can show in the same way as in
\cite[Theorem 3.4]{Cie2} that
\begin{equation}\label{eq:orliczineqL}
\|f\|_{\mathfrak{L}^{\Phi_{\beta}}(I;X)} \leq
\|f\|_{L^{\Phi_{\beta}}(I;X)}.
\end{equation}


\subsection{Gaussian random variables}
Let $(\O, \mathcal A, \P)$ be a complete probability space. As in
\cite{LeTa} let $X$ be Banach space with the following property:
there exists a sequence $(x_n^*)_{n\geq 1}$ in $X^*$ such that
$\|x_n^*\|\leq 1$ and $\|x\| = \sup_{n\geq 1} |x_n^*(x)|$. Such a
Banach space will be said to {\em admit a norming sequence of
functionals}. Examples of such Banach spaces are all separable
Banach spaces, but also spaces like $l^\infty$. As in \cite{LeTa} a
mapping $\xi:\O\to X$ will be called {\em a centered Gaussian} if
for all $x^*\in \text{span}\{x_n^*:n\geq 1\}$ the random variable
$\lb \xi, x^*\rb$ is a centered Gaussian. For a centered Gaussian
random variable we define
\begin{equation}\label{eq:sigmaUnitBall}
  \sigma(\xi) = \sup_{n\geq 1} (\E|\lb \xi, x^*_n\rb|^2)^{1/2}.
\end{equation}
In \cite{LeTa} it is proved that
\[
  \lim_{t\to \infty} \frac{1}{t^2}\log \P(\|X\|>t) = -\frac{1}{2\sigma^2},
\]
so that the value of $\sigma$ is independent of the norming sequence
$(x_n^*)_{n\geq 1}$.

We make some comment on the above definition of a Gaussian random
variable. We do not assume that $\xi$ is a Borel measurable mapping.
The only obvious fact we will use is that the mapping $\omega\mapsto
\|\xi(\omega)\|$ is measurable. If $\xi$ is a Gaussian random
variable that takes values in a separable subspace of $X$, then
$\xi$ is Borel measurable and one already has that $\lb \xi, x^*\rb$
is a centered Gaussian random variables for all $x^*\in X^*$.

A random variable $\xi:\O\to X$ is called {\em tight} if the measure
$\P\circ \xi^{-1}$ is tight, and it is called {\em Radon} if
$\P\circ \xi^{-1}$ is Radon. If $X$ is a separable Banach space,
then every Borel measurable random variable $\xi:\O\to E$ is Radon,
and in particular tight. Conversely, if a Gaussian random variable
$\xi:\O\to X$ is tight, then it almost surely takes values in a
separable subspace of $X$. The next result is well-known and a short
proof can be found in \cite[p.\ 61]{LeTa}.

\begin{proposition}\label{prop:nondeg}
Let $X$ be a Banach space and let $\xi:\O\to X$ be a centered
Gaussian. If $\xi$ is tight, then $\P(\|\xi\|<r)>0$ for all $r>0$.
\end{proposition}



%
%

\section{Maximal estimates for sequences of Gaussian random variables}\label{sec:linfty}


%
%
%
%

The next proposition together with Theorem \ref{thm:total} may be
considered as the vector-valued extension of a result in \cite{Hy}.

\begin{proposition}\label{prop:supestX}
Let $X$ be a Banach space which admits a norming sequence of
functionals $(x_n^*)_{n\geq 1}$. Let $\Theta$ be as in
\eqref{Phiexp}. Let $(\xi_n)_{n\geq 1}$ be $X$-valued centered
Gaussian random variables with first moments and weak variances
\[m_n = \E\|\xi_n\|,\]
\[\sigma_n = \sup_{m\geq 1} (\E|\lb \xi_n, x^*_m\rb|^2)^{1/2}.\]
It holds that
\[
\E\sup_{n\geq 1} \|\xi_n\| \leq m + 3
\rho_{\Theta}((\sigma_n)_{n\geq 1}),
\]
where $m=\sup_{n\geq 1}m_n$.

Moreover, if any linear combination of the $(\xi_n)_{n\geq 1}$ is a
Gaussian random variable and if $\E\sup_{n\geq 1} \|\xi_n\|<\infty$,
then $\xi:=(\xi_n)_{n\geq 1}$ is an $l^\infty(X)$-valued Gaussian
random variable.
\end{proposition}

By the Kahane-Khinchine inequalities (cf. \cite[Corollary
3.4.1]{KwWo}) one obtains a similar estimate for the $p$-th moments
of $\sup_{n\geq 1} \|\xi_n\|$. However, this also follows by
extending the proof below.

\begin{proof}
We may write
\[\E\sup_{n\geq 1} \|\xi_n\| \leq \E\sup_{n\geq 1} |\|\xi_n\| - m_n| + \sup_{n\geq 1} m_n \]
By \cite[(3.2)]{LeTa} for all $t>0$, we have
\begin{equation}\label{eq:LeTa32}
  \P(|\|\xi_n\| - m_n|>t)\leq 2\exp\Big(-\frac{t^2}{2\sigma_n^2}\Big).
\end{equation}
For each $\delta>0$ it follows that {\small
\begin{equation}\begin{aligned}
  \E \sup_{n\geq 1} & |\|\xi_n\| - m_n|
  = \int\limits_0^\infty \P(\sup_{n\geq 1}|\|\xi_n\| - m_n|>t) \, dt
 \\ & \leq \delta + \int\limits_\delta^\infty \P(\sup_{n\geq 1} |\|\xi_n\| - m_n|>t) \, dt
  \leq \delta + \sum_{n\geq 1} \int\limits_\delta^\infty \P(|\|\xi_n\| - m_n|>t) \, dt \\
 & \leq \delta + \sum_{n\geq 1} 2\int\limits_\delta^\infty  \exp\Big(-\frac{t^2}{2\sigma_n^2}\Big) \, dt
  = \delta + \sum_{n\geq 1} 2\int\limits_{\delta/\sigma_n}^\infty  \sigma_n\exp\Big(-\frac{t^2}{2}\Big) \, dt \\
 & \leq \delta + 2\sum_{n\geq 1} \frac{\sigma_n^2}{\delta} \exp\Big(-\frac{\delta^2}{2\sigma_n^2}\Big)
  = \delta\Big[1 + 2\sum_{n\geq 1}\frac{\sigma_n^2}{\delta^2}
    \exp\Big(-\frac{\delta^2}{2\sigma_n^2}\Big)\Big]
\end{aligned}\end{equation}}
where we used the standard estimate
\[\int\limits_\delta^\infty e^{-t^2/2} \, dt \leq \frac{1}{\delta}\exp(-\delta^2/2).\]

If $\delta>0$ is chosen so that the last series sums up to at most
$1$, then we have shown that $\E\sup_{n\geq 1} |\|\xi_n\| - m_n|\leq
3\delta$. Taking the infimum over all such $\delta$, we obtain the
result.

The final assertion follows from the definition of a Gaussian random
variable using the norming sequence of functionals $(e_m \otimes
x^*_n)_{m,n\geq 1}$.
\end{proof}

\begin{remark}\label{rem:expVsLp}
The infimum appearing in Proposition~\ref{prop:supestX} is dominated
by
\[
  \Big[\Big(\frac{p-1}{e}\Big)^{\frac{p-1}{2}}\sum_{n\geq 1}\sigma_n^{p+1}\Big]^{1/(p+1)}
\]
for any $p\in\left[1,\infty\right[$. (Interpret $0^0=1$ for $p=1$.)
This follows from the elementary estimate
$e^{-x^2/2}\leq[(p-1)/e]^{(p-1)/2}x^{1-p}$ applied to
$x=\delta/\sigma_n$.
\end{remark}

For an $X$-valued random variable $\xi$ we take a median $M$ such
that
\[\P(\|\xi\|\leq M)\geq 1/2 \ \ \ \text{and} \ \ \ \P(\|\xi\|\geq M)\geq 1/2.\]
For convenience we will take $M=M(\xi)$ to be the smallest possible
$M$. Notice that for all $p\in (0, \infty)$, $\E\|\xi\|^p\geq
\frac{M^p}{2}$.

Alternatively, we could have replaced the estimate~\eqref{eq:LeTa32}
in the above proof by
\[
  \P(|\|\xi\|-M|>t)\leq \exp\Big(-\frac{t^2}{2\sigma^2}\Big).
\]
(see \cite[Lemma 3.1]{LeTa}) to obtain:

\begin{proposition}\label{prop:X}
Let $X$ be a Banach space which admits a norming sequence of
functionals $(x_n^*)_{n\geq 1}$. Let $\Theta$ be as in
\eqref{Phiexp}. Let $(\xi_n)_{n\geq 1}$ be $X$-valued centered
Gaussian random variables with medians $M_n$ and
weak variances
\[\sigma_n = \sup_{m\geq 1} (\E|\lb \xi_n, x^*_m\rb|^2)^{1/2}.\]
It holds that
\[\E\sup_{n\geq 1} \|\xi_n\| \leq M + 2\rho_{\Theta}((\sigma_n)_{n\geq 1}),\]
where $M=\sup_{n\geq 1} M_n$.
\end{proposition}

If the $\xi_n$ are independent Gaussian random variables, then a
converse to Proposition \ref{prop:supestX} holds.

\begin{theorem}\label{thm:total}
Let $X$ be a Banach space which admits a norming sequence of
functionals. Let $\Theta$ be as in \eqref{Phiexp}. Let
$(\xi_n)_{n\geq 1}$ be $X$-valued independent centered Gaussian
random variables with first moments $(m_n)_{n\geq 1}$ and weak
variances $(\sigma_n)_{n\geq 1}$. Let $m=\sup_{n\geq 1}m_n$. It
holds that
\[\begin{aligned}
\E\sup_{n\geq 1} \|\xi_n\| &\eqsim m+\rho_{\Theta}((\sigma_n)_{n\geq
1})
\\ & \eqsim m+\|(\sigma_n)_{n\geq 1}\|_{\Theta}.
\end{aligned}\]
Moreover, if one of these expressions is finite, then
$\xi:=(\xi_n)_{n\geq 1}$ is an $l^\infty(X)$-valued Gaussian random
variable.
\end{theorem}

Recall from Subsection \ref{subsec:Orlicz} and the definition of
$\Theta$ that
\[\|(\sigma_n)_{n\geq 1}\|_{\Theta} = \inf_{\delta>0}\Big\{\frac{1}{\delta}\Big[1+ \sum_{n\geq 1} \delta^2\sigma_n^2
\exp\Big(-\frac{1}{2\delta^2\sigma_n^2}\Big)\Big] \Big\}.
\]

\begin{proof}
The second two sided estimate follows from Lemma
\ref{lem:OrliczLux}.

The estimate $\lesssim$ in the first comparison has been obtained in
Proposition \ref{prop:supestX}. To prove $\gtrsim$, note that
$\E\sup_{n\geq 1} \|\xi_n\|\geq m$ is clear. As for the estimate for
$\rho_{\Theta}((\sigma_n)_{n\geq 1})$, by scaling we may assume that
$\E\sup_{n\geq 1} \|\xi_n\|$=1. Then one has $\P(\sup_{n\geq 1}
\|\xi_n\|>3)\leq 1/3$, and therefore
\[\begin{aligned}
1/3\leq \P(\sup_{n\geq 1} \|\xi_n\|\leq 3) &= \prod_{n\geq 1}
\P(\|\xi_n\|\leq 3) = \prod_{n\geq 1} (1-\P(\|\xi_n\|> 3))\\ &\leq
\prod_{n\geq 1} \exp\Big(-\P(\|\xi_n\|> 3)\Big).
\end{aligned}\] It follows that
\[\log 3\geq \sum_{n\geq 1}\P(\|\xi_n\|> 3).\]
Let $\e\in (0,1)$ be an arbitrary number. If for each $n\geq 1$, we
choose $k_n$ such that $(\E\lb \xi_n, x^*_{k_n}\rb^2)^{1/2} \geq
\sigma_n(1-\e)$, then we obtain
\[\begin{aligned}
  \log 3 &\geq \sum_{n\geq 1}\P(\|\xi_n\|> 3)
  \geq \sum_{n\geq 1}\P(|\lb \xi_n,x_{k_n}^*\rb|> 3) \\
  &\geq \sqrt{\frac{2}{\pi}} \sum_{n\geq 1}\frac{3\sigma_n(1-\e)}{\sigma_n^2(1-\e)^2+9}
    \exp\Big(-\frac{9}{2\sigma_n^2(1-\e)^2}\Big)
\end{aligned}\]
where we used $\int\limits_a^\infty e^{-t^2/2} \, dt\geq
\frac{a}{1+a^2}e^{-a^2/2}$. Next, we have
\[
  \sigma_n^2 = \sup_{m\geq 1} \E\lb \xi_n, x_m^*\rb^2
  =\frac{\pi}{2}\sup_{m\geq 1} \E|\lb \xi_n, x_m^*\rb|\leq \frac{\pi}{2}\E\|\xi_n\|\leq \frac{\pi}{2},
\]
hence $\sigma_n^2(1-\e)^2+9\leq\pi/2+9<11$ and
$\sqrt{2/\pi}\cdot\sigma_n\geq 2/\pi\cdot\sigma_n^2$, thus
\[
  \log 3 \geq \frac{6}{11\pi}\sum_{n\geq 1}\sigma_n^2(1-\e) \exp\Big(-\frac{9}{2\sigma_n^2(1-\e)^2}\Big).
\]
This being true for all $\e>0$, it follows in the limit that
\[
  \sum_{n\geq 1}\Big(\frac{\sigma_n}{3}\Big)^2 \exp\Big(-\frac{9}{2\sigma_n^2}\Big)
  \leq\log 3\cdot \frac{11\pi}{6\cdot 9}<1.
\]
Therefore, $\rho_\Theta((\sigma_n)_{n\geq 1})\leq 3$.

The last assertion follows as in Proposition \ref{prop:supestX}.
\end{proof}

From the proof of Theorem \ref{thm:total} we actually see that
\[\E\sup_{n\geq 1} \|\xi_n\| \geq \max\Big\{\frac{1}{3}\rho_{\Theta}((\sigma_n)_{n\geq 1}), m\Big\}.\]


\begin{remark}
A similar proof as presented above shows that the function $\Theta$
in Theorem \ref{thm:total} can be replaced by the (non-convex)
function $\Phi$ defined below Lemma \ref{lem:OrliczLux}. Since we
prefer to have an Orlicz space, we use the convex function $\Theta$.
\end{remark}

In the real-valued case, $m$ is not needed in the estimate of
Theorem \ref{thm:total}. This is due to the fact that it can be
estimated by $\sup_{n\geq 1}\sigma_n$. The following simple example
shows that in the infinite dimensional setting this is not the case.
We shall also encounter the same phenomenon in a more serious
example in the proof of Theorem~\ref{thm:momentest}.

\begin{example}\label{ex:mVsSigma}
Let $p\in [1, \infty]$ and let $X=l^p$ with the standard unit
vectors denoted by $e_n$. Let $(\sigma_n)_{n\geq 1}$ be a sequence
of positive real numbers with
\begin{equation*}
  m_p:=\Big(\sum_{n\geq 1} \sigma_n^p\Big)^{1/p} <\infty \ \ \text{if} \ \ p<\infty
\end{equation*}
and
\begin{equation*}
  m_{\infty}:=\rho_{\Theta}((\sigma_n)_{n\geq 1})<\infty \ \
  \text{if} \ \ p=\infty.
\end{equation*}
Let $(\g_n)_{n\geq 1}$ be a sequence of independent standard
Gaussian random variables. Then $\xi = \sum_{n\geq 1} \sigma_n \g_n
e_n$ defines an $X$-valued Gaussian random variable with \(m(\xi) =
\E\|\xi\|\eqsim_p m_p\) and
\[\sigma(\xi) =
\left\{\begin{array}{cc}
  \sup_{n\geq 1} \sigma_n & p\in [2, \infty] \\
  \Big(\sum_{n\geq 1} \sigma_n^{r}\Big)^{\frac1r} & p\in [1, 2),
\end{array}
\right.
\]
where $r= \frac{2p}{2-p}$.
\end{example}




\section{Besov regularity of Brownian paths}\label{sec:regBM}

We say that an $X$-valued process $(W(t))_{t\in [0,1]}$ is a {\em
Brownian motion} if it is strongly measurable and for all $x^*\in
E^*$, $(\lb W(t), x^*\rb)_{t\in [0,1]}$ is a real Brownian motion
starting at zero. Let $Q$ be the covariance of $W(1)$. For the
process $W$ we have
\begin{enumerate}
\item $W(0) = 0$,
\item $W$ has a version with continuous paths,
\item $W$ has independent increments,
\item For all $0\leq s<t<\infty$, $W(t) - W(s)$ has distribution
$\mathcal{N}(0,(t-s)Q)$.
\end{enumerate}
In this situation we say that $W$ is a Brownian motion {\em with
covariance $Q$}. Notice that every process $W$ that satisfies (3)
and (4) has a path-wise continuous version (cf.\ \cite[Theorem
3.23]{Kal}).

In the next result we obtain a Besov regularity result for Brownian
motions. The case of real valued Brownian motions has been
considered in \cite{Cie,Cie2,Roy}. But even in the real-valued case
we believe the proof is new and more direct.

\begin{theorem}\label{thm:pathsasB}
Let $X$ be a Banach space and let $p,q\in [1, \infty)$. For an
$X$-valued non-zero Brownian motion $W$ we have
\begin{equation*}\begin{split}
  W\in B^{1/2}_{\Phi_2,\infty}(0,1;X)\subset B^{1/2}_{p,\infty}(0,1;X)\qquad a.s.,\\
  W\notin B^{1/2}_{p,q}(0,1;X)\qquad a.s.
\end{split}\end{equation*}
\end{theorem}

\begin{proof}
Denote
\[Y_{n,p}:=2^{n/2}\|W(\cdot+2^{-n})-W\|_{L^p(I(2^{-n});X)}.\]
We may write
\begin{equation*}\begin{aligned}
Y^p_{n,p} &= \int\limits_0^{1-2^{-n}} 2^{np/2}\|W(t+2^{-n}) -
W(t)\|^p \, dt
\\ & = \sum_{m=1}^{2^n-1} \int\limits_{(m-1) 2^{-n}}^{m 2^{-n}} 2^{np/2}\|W(t+2^{-n}) - W(t)\|^p \, dt
\\ & = \sum_{m=1}^{2^n-1} 2^{-n}\int\limits_{0}^{1} 2^{np/2}\|W((s+m)2^{-n}) - W((s+m-1)2^{-n})\|^p \,
ds
\\ & =\int\limits_{0}^{1}  2^{-n}\sum_{m=1}^{2^n-1} \|\g_{n,m,s}\|^p \, ds
\end{aligned}\end{equation*}
Here \(\g_{n,m,s} = 2^{n/2} (W((s+m)2^{-n}) - W((s+m-1)2^{-n}))\).
For fixed $s\in (0,1)$ and $n\geq 1$, $(\g_{n,m,s})_{m\geq 1}$ is a
sequence of independent random variables distributed as $W(1)$.
Denote $c_p = (\E\|W(1)\|^p)^{1/p}$. If we take second moments we
may use Jensen's inequality to obtain
\begin{equation*}\begin{aligned}
\E \Big( Y^p_{n,p} - c_p^p\Big)^2 &= \E \Big|\int\limits_{0}^{1}
\big[2^{-n} \sum_{m=1}^{2^n-1} \|\g_{n,m,s}\|^p - c_p^p\big]
\,ds\Big|^2
\\ & \leq \int\limits_{0}^{1}  \E \Big| 2^{-n} \sum_{m=1}^{2^n-1}\big(\|\g_{n,m,s}\|^p -c_p^p\big)-2^{-n}c_p^p\Big|^2 \,ds
\\ & = \int\limits_{0}^{1}\Big[  2^{-2n}(2^n-1) (c_{2p}^{2p} -c_p^{2p})+2^{-2n}c_p^{2p} \Big] ds
\\ &= 2^{-n}\big[(1-2^{-n})c_{2p}^{2p}-(1-2^{1-n})c_p^{2p}\big].
\end{aligned}\end{equation*}
It follows that for a fixed $\e>0$, we have
\begin{equation*}\begin{aligned}
  \sum_{n\geq 1} &\P\Big( \big|Y^p_{n,p} - c_p^p\big| >\e\Big) \leq
  \frac{1}{\e^2}\sum_{n\geq 1}\E \Big( Y^p_{n,p} - c_p^p\Big)^2<\infty,
\end{aligned}\end{equation*}
which implies, by the Borel--Cantelli Lemma, that
\begin{equation*}
  \P\Big(\big|Y^p_{n,p} - c_p^p\big| >\e \
  \text{infinitely often}\Big) = 0.
\end{equation*}
This in turn gives that
\begin{equation}\label{eq:limBesovW}
  \lim_{n\to \infty} 2^{n/2}\|W(\cdot+2^{-n}) - W\|_{L^p(I(2^{-n});X)}
   = (\E\|W(1)\|^p)^{1/p} \ \ \text{a.s.}
\end{equation}
This shows immediately that the paths are a.s.\ in
$B^{1/2}_{p,\infty}(0,1;X)$. From the above calculation it is also
clear that $W\notin B^{1/2}_{p,q}(0,1;X)$ a.s.\ for $q\in [1,
\infty)$. Next we show that the paths are in
$B^{1/2}_{\Phi_2,\infty}(0,1;X)$ a.s. Note that
$(\E\|W(1)\|^p)^{1/p}\eqsim p^{1/2}$ as $p\to\infty$. The upper
estimate $\lesssim$ is a consequence of Fernique's theorem (which
says that $\|W(1)\|^2$ is exponentially integrable, since $W(1)$ is
a non-zero $X$-valued Gaussian random variable), whereas $\gtrsim$
follows from the corresponding estimate for real Gaussians after
applying a functional. We proved that $\E(Y^p_{n,p}-c_p^p)^2\leq
c_{2p}^{2p} 2^{-n}$. Therefore,
\[
  \E(Y^p_{n,p} c_p^{-p}-1)^2\leq 2^{-n}c_{2p}^{2p} c_p^{-2p}
  \leq 2^{-n} K^{2p},
\]
where $K\geq 1$ is some constant. Hence for all $\lambda>1$,
\[  \P(Y_{n,p} c_p^{-1}>\lambda)
  \leq
  \P(|Y^p_{n,p} c_p^{-p}-1|>\lambda^p-1)
  \leq 2^{-n}K^{2p} (\lambda^p-1)^{-2},
\]
and thus for $\lambda = 2K$
\[
  \sum_{n,p=1}^{\infty}\P(Y_{n,p} c_p^{-1}>\lambda)
  \leq \sum_{n=1}^{\infty}2^{-n}\sum_{p=1}^{\infty}K^{2p} (\lambda^p-1)^{-2}
  <\infty
\]
so that by the Borel--Cantelli lemma
\[
  \P\Big(Y_{n,p} c_p^{-1}>\lambda\text{ for infitely many pairs }(n,p)\Big)=0.
\]
Since $c_p\eqsim p^{1/2}$ this means that a.s.
\[
  \sup_{n,p} 2^{n/2}\|W(\cdot+2^{-n})-W\|_{L^p(I(2^{-n});X)}p^{-1/2}<\infty.
\]
\end{proof}

\section{Brownian motions as random variables in Besov spaces}\label{sec:rvW}

From the pathwise properties of $W$ studied in the previous section,
we know that we have a function $W:\Omega\to B^{1/2}_{p,\infty}$. We
now go into the measurability issues in order to promote it to a
random variable.

\begin{theorem}\label{thm:WGaussian}
Let $X$ be a Banach space and let $p\in [1, \infty)$. Then an
$X$-valued Brownian motion $W$ is a
$B^{1/2}_{p,\infty}(0,1;X)$-valued, and even
$B^{1/2}_{\Phi_2,\infty}(0,1;X)$-valued, Gaussian random variable.
In particular, there exists an $\e>0$ such that
\begin{equation*}
  \E\exp\big(\e\|W\|_{B^{1/2}_{\Phi_2,\infty}(0,1;X)}^2\big)<\infty.
\end{equation*}
If $W$ is non-zero, then the random variables $W:\O\to
B^{1/2}_{p,\infty}(0,1;X)$ and $W:\Omega\to
B^{1/2}_{\Phi_2,\infty}(0,1;X)$ are not tight. In fact,
\begin{equation*}
  \tau_1 := \inf\{\lambda\geq 0: \P(\|W\|_{B^{1/2}_{p,\infty}(0,1;X)}\leq \lambda)>0\} \geq (\E\|W(1)\|^p)^{1/p},
\end{equation*}
and consequently also
\begin{equation*}
  \tau_2:=\inf\{\lambda\geq 0: \P(\|W\|_{B^{1/2}_{\Phi_2,\infty}(0,1;X)}\leq \lambda)>0\}>0.
\end{equation*}
\end{theorem}

There is some interest in the numbers $\tau_1$ and $\tau_2$. For
general theory we refer the reader to \cite[Chapter 3]{LeTa}.

For the proof we need the following easy lemma.

\begin{lemma}\label{lem:Besovnormingsq}
Let $X$ be a Banach space which admits a norming sequence, let
$0<\alpha<1$ and $0<\beta<\infty$. Then for all $p\in[1,\infty)$
there exist
\begin{equation*}
  (\Lambda_{pjk})_{j\geq 0,k\geq 1}
  \subset B^{\alpha}_{p,\infty}(0,1;X)^*
  \subset B^{\alpha}_{\Phi_{\beta},\infty}(0,1;X)^*,
\end{equation*}\begin{equation*}
  (f_{pjk})_{j\geq 0,k\geq 1}
  \subset C^{\infty}([0,1];X^*),
\end{equation*}
such that: for all $\phi\in B^{\alpha}_{p,\infty}(0,1;X)$ there are
the representations
\begin{equation*}
  \lb \phi, \Lambda_{p0k}\rb = \int\limits_0^1 \lb \phi(t), f_{p0k}(t)\rb \, dt, \ \ k\geq 1,
\end{equation*}
\begin{equation*}
  \lb \phi, \Lambda_{pjk}\rb = \int\limits_0^{1-2^{-j}}  2^{j\alpha}\lb \phi(t+2^{-j}) - \phi(t), f_{pjk}(t)\rb \, dt,\qquad j,k\geq 1;
\end{equation*}
we have the upper norm bounds
\begin{equation*}
  p^{-1/\beta}\|\Lambda_{pjk}\|_{B^{\alpha}_{\Phi_{\beta},\infty}(0,1;X)^*}
  \leq \|\Lambda_{pjk}\|_{B^{\alpha}_{p,\infty}(0,1;X)^*}\leq 1,\qquad k\geq 1;
\end{equation*}
and finally the sequences are norming in the following sense:
\begin{equation*}\begin{split}
  \|\phi\|_{B^{\alpha}_{p,\infty}(0,1;X)} &= \sup_{j\geq 0,k\geq 1} |\lb\phi,\Lambda_{pjk}\rb|, \\
  \|\phi\|_{B^{\alpha}_{\Phi_{\beta},\infty}(0,1;X)}
   &= \sup_{p\geq 1,j\geq 0,k\geq 1}p^{-1/\beta} |\lb\phi,\Lambda_{pjk}\rb|.
\end{split}\end{equation*}
\end{lemma}

\begin{proof}
Let $(x_n^*)_{n\geq 1}$ be a norming sequence for $X$. Let
$I=[a,b]$. First observe that there exists a sequence
$(F_{k})_{k\geq 1}$ in $L^{p'}(I;X^*)$, with norm smaller than or
equal to one, which is norming for $L^{p}(I;X)$. Such a sequence is
easily constructed using the $(x_n^*)_{n\geq 1}$ and standard
duality arguments. By an approximation argument we can even take the
$(F_k)_{k\geq 1}$ in $C^{\infty}(I;X^*)$.

To prove the lemma, let first $a=0$ and $b=1$, and let
$(f_{p0k})_{k\geq 1}$ be the above constructed sequence
$(F_k)_{k\geq 1}$. Next we fix $j\geq 1$ and let $a=0$ and
$b=1-2^{-j+1}$ and let $(f_{pjk})_{k\geq 1}$ be the above
constructed sequence for this interval. Let $\Lambda_{pjk}$ be the
elements in $B^{\alpha}_{p,\infty}(0,1;X)^*$ defined as in the
statement in the lemma. It is easily checked that this sequence
satisfies the required properties.
\end{proof}

\begin{proof}{Proof of Theorem {\ref{thm:WGaussian}}}
Since $W$ is strongly measurable as an $X$-valued process we may
assume that $X$ is separable and therefore that it admits a norming
sequence. In Theorem \ref{thm:pathsasB} it has been shown that the
paths of $W$ are a.s.\ in $B^{1/2}_{\Phi_2,\infty}(0,1;X)\subset
B^{1/2}_{p,\infty}(0,1;X)$ for all $p\in[1,\infty)$. It follows from
Lemma \ref{lem:Besovnormingsq} that there exists a norming sequence
of functionals $(\Lambda_n)_{n\geq 1}$ for
$B^{1/2}_{\Phi_2,\infty}(0,1;X)$, as well as in each
$B^{1/2}_{p,\infty}(0,1;X)$, such that $\lb W, \Lambda\rb$ is a
centered Gaussian random variable for all $\Lambda\in \text{span}
\{\Lambda_n, n\geq 1\}$. Therefore, by definition it follows that
$W$ is a centered Gaussian random variable. The exponential
integrability follows from \cite[Corollary 3.2]{LeTa}.

The last assertion follows from~\eqref{eq:limBesovW}.
This also shows that $W$ is not tight since, by Proposition \ref{prop:nondeg}, for centered Gaussian
measures which are tight, one has $\tau=0$.
\end{proof}

\section{Moment estimates for Brownian motions in Besov spaces}\label{sec:moment}

Now that we know that
\[
  \E\|W\|_{B^{1/2}_{p,\infty}(0,1;X)}<\infty,\qquad
  \E\|W\|_{B^{1/2}_{\Phi_2,\infty}(0,1;X)}<\infty,
\]
it seems interesting to estimate these quantities. For this we need
a convenient representation of $X$-valued Brownian motions.

Recall that a family $W_H=(W_H(t))_{t\in \R_+}$ of bounded linear
operators from $H$ to $L^2(\O)$ is called an {\em $H$-cylindrical
Brownian motion} if
\begin{enumerate}
\item $W_H h = (W_H(t)h)_{t\in \R_+}$ is a real-valued Brownian motion for each $h\in H$,
\item $ \E (W_H(s)g \cdot W_H(t)h) = (s\wedge t)\,[g,h]_{H}$ for all $s,t\in \R_+, \ g,h\in H.$
\end{enumerate}
We always assume that the $H$-cylindrical Brownian motion $W_H$ is
adapted to a given filtration $\F$, i.e., the Brownian motions
$W_Hh$ are adapted to $\F$ for all $h\in H$. Notice that if
$(h_n)_{n\geq 1}$ is an orthonormal basis for $H$, then $(W_H
h_n)_{n\geq 1}$ are independent standard real-valued Brownian
motions.

Let $W:\R_+\times\O\to E$ be an $E$-valued Brownian motion and let
$Q\in\calL(E\s,E)$ be its covariance operator. Let $H_Q$ be the
reproducing kernel Hilbert space or Cameron--Martin space (cf.\
\cite{Bog,VTC}) associated with $Q$ and let $i_W: H_Q\embed E$ be
the inclusion operator. Then the mappings
\[
 W_{H_Q}(t):  i_W\s x\s \mapsto \lb W(t),x\s\rb
\]
uniquely extend to an $H_Q$-cylindrical Brownian motion $W_{H_Q}$,
so that in particular
\begin{equation}\label{eq:repW}
  \lb W(t),x^*\rb = W_{H_Q}(t)i_W^*x^*.
\end{equation}

\begin{lemma}\label{lem:extra}
There holds, for all $p\in[1,\infty)$,
\[
  \|i_W\|=\sigma(W(1))\lesssim\frac{1}{\sqrt{p}}(\E\|W(1)\|^p)^{1/p}.
\]
\end{lemma}

\begin{proof}
Note first that, since $\lb W(t),x^*\rb$ is a real-valued Gaussian
random variable, its moments satisfy
\begin{equation}\label{eq:realGauss}
  (\E\|\lb W(t),x^*\rb|^p)^{1/p}=\gamma_p(\E|\lb W(t),x^*\rb|^2)^{1/2},
\end{equation}
where the $\gamma_p$ are universal constants behaving like
$\gamma_p\eqsim\sqrt{p}$ for $p\in[1,\infty)$.

On the other hand, by \eqref{eq:repW} and the definition of
cylindrical Brownian motion,
\[
  (\E|\lb W(t),x^*\rb|^2)^{1/2}=\sqrt{t}\|i_W^*x^*\|.
\]
With $t=1$, taking supremum over all $x^*\in X^*$ of unit norm, and
recalling that $\|i_W\|=\|i_W^*\|$, this proves the first equality
in the assertion. The second then follows from \eqref{eq:realGauss}
and the obvious estimate
\[
  (\E\|\lb W(t),x^*\rb|^p)^{1/p}\leq(\E\|W(t)\|^p)^{1/p}
\]
for $\|x^*\|\leq 1$.
\end{proof}

\begin{lemma}\label{lem:varianceW}
Let $c>0$, and $J\subset\R_+$ be an interval of length $|J|\geq c$.
Consider $W(\cdot+c)-W$ as an $L^p(J,X)$-valued Gaussian random
variable. Then
\[
  \sigma(W(\cdot+c)-W)\eqsim c^{1/2+1/p}\|i_W\|.
\]
\end{lemma}


\begin{proof}
To prove the claim take $f\in L^{p'}(J;X^*)$. We also use the same
symbol for its extension to $\R$ with zero fill. The representation
\eqref{eq:repW}, the Stochastic Fubini theorem, and the It\^o
isometry yield
\begin{equation*}\begin{aligned}
  \Big(\E\Big|\int_J &\lb (W(t+c) - W(t)),f(t)\rb \,dt \Big|^2\Big)^{1/2} \\
  & = \Big(\E\Big|\int_J (W_H(t+c) - W_H(t)) i_W^*f(t) \,dt \Big|^2\Big)^{1/2} \\
  & = \Big(\E\Big|\int_{\R} \int_{\R_+} \one_{[t,t+c]}(s) i_W^* f(t) \, d W_H(s)\, dt \Big|^2\Big)^{1/2} \\
  & = \Big(\E\Big|\int_{\R_+} \one_{[0,c]}* (i_W^* f)(s) \, dW_H(s) \Big|^2\Big)^{1/2} \\
  & = \Big(\int_{\R}\|\one_{[0,c]}* (i_W^* f)(s)\|_H^2 \, ds\Big)^{1/2}.
\end{aligned}\end{equation*}
Taking the supremum over all $f\in L^{p'}(J;X^*)$ of unit norm, we
find that
\[
  \sigma(W(\cdot+c)-W)=\|(\one_{[0,c]}*)\otimes i_W^*\|_{L^{p'}(J;X^*)\to L^2(\R;H)}.
\]
By Young's inequality with $1+1/2=1/p'+1/r$ it follows that the
operator norm is dominated by
\[
  \|\one_{[0,c]}\|_{L^r}\|i_W^*\|_{X^*\to H}
  =c^{1/p+1/2}\|i_W\|.
\]
On the other hand, if we test with the functions $f=\one_I\otimes
x^*\in L^{p'}(J;X^*)$, where $I\subseteq J$ has length $c$, we find
that
\[\begin{split}
  \|\one_{[0,c]}* (i_W^* f)\|_{L^2(H)}
  &=\|\one_{[0,c]}*\one_I\|_{L^2}\|i_W^* x^*\|_H \\
  &=(2/3)^{1/2} c^{3/2}\|i_W^* x^*\|_H
  \eqsim c^{1/2+1/p}\frac{\|i_W^*x^*\|_H}{\|x^*\|_{X^*}} \|f\|_{L^{p'}(X^*)}.
\end{split}\]
Taking the supremum over $x^*\in X^*\setminus\{0\}$ we get the other
side of the asserted norm equivalence.
\end{proof}

\begin{corollary}\label{cor:varianceW}
Let $c\in(0,e^{-1/2}]$, and $J\subset\R_+$ be an interval of length
$|J|\geq c$. Consider $W(\cdot+c)-W$ as an
$\mathfrak{L}^{\Phi_2}(J;X)$-valued Gaussian random variable. Then
\[
  \sigma(W(\cdot+c)-W)\eqsim (\log c^{-1})^{-1/2}c^{1/2}\|i_W\|.
\]
\end{corollary}

\begin{proof}
We note that the functionals $p^{-1/2}\Lambda_{p0k}$ from
Lemma~\ref{lem:Besovnormingsq} (with $\beta=2$) provide a norming
sequence for $\mathfrak{L}^{\Phi_2}(0,1;X)$, and the same
construction can be adapted to another interval. Hence
\[\begin{split}
  \sigma_{\mathfrak{L}^{\Phi_2}(J;X)}&(W(\cdot+c)-W) \\
  &=\sup_{p\geq 1}p^{-1/2}\sup_{k\geq 1}\Big(\E\Big|\int_J \lb (W(t+c)-W(t)), f_{p0k}(t)\rb \, dt\Big|^2\Big)^{1/2}\\
  &=\sup_{p\geq 1}p^{-1/2}\sigma_{L^p(J;X)}(W(\cdot+c)-W)\\
  &\eqsim\sup_{p\geq 1}p^{-1/2}c^{1/2+1/p}\|i_W\|\\
  &\eqsim(\log c^{-1})^{-1/2}c^{1/2}\|i_W\|,
\end{split}\]
where an elementary maximum value problem was solved in the last
step.
\end{proof}

\begin{theorem}\label{thm:momentest}
Let $X$ be a Banach space. Let $p\in [1, \infty)$. For an $X$-valued 
Brownian motion $W$ we have
\begin{equation}\label{eq:expBesov}
  \E\|W\|_{B^{1/2}_{p,\infty}(0,1;X)}
  \eqsim  (\E\|W(1)\|^p)^{1/p},
\end{equation}
\begin{equation}\label{eq:expBesovOrlicz}
  \E\|W\|_{B^{1/2}_{\Phi_2,\infty}(0,1;X)}
  \eqsim \E\|W(1)\|.
\end{equation}
\end{theorem}

\begin{remark}
By \cite[Corollary~3.2]{LeTa}, the estimate \eqref{eq:expBesov}
implies that
\[
   \E\|W\|_{B^{1/2}_{p,\infty}(0,1;X)}\lesssim\sqrt{p}\,\E\|W(1)\|,
\]
but we do not know if there is a two sided comparison here. The
above estimate is also an immediate consequence of
\eqref{eq:expBesovOrlicz} and the definition of the various norms.
\end{remark}

\begin{proof}
As in Theorem \ref{thm:WGaussian} we may assume that $X$ admits a
norming sequence.

The estimate $\gtrsim$ in \eqref{eq:expBesov} follows from
\eqref{eq:limBesovW}. Let us then consider the other direction.
Clearly,
\[
  \E\|W\|_{L^p(0,1;X)}
  \leq (\E\|W\|_{L^{\infty}(0,1;X)}^2)^{1/2}
  \leq 2(\E\|W(1)\|^2)^{1/2}
  \lesssim \E\|W(1)\|
\]
by Doob's maximal inequality and the equivalence of Gaussian
moments. Next we consider
\begin{equation}\label{eq:supW}
  \E\sup_{j\geq 1} 2^{j/2} \|W(\cdot+2^{-j}) -
  W\|_{L^p(0,1-2^{-j};X)}.
\end{equation}
This can be estimated using Proposition \ref{prop:supestX} with the
$L^p(0,1;X)$-valued Gaussian random variables
$\xi_j=2^{j/2}[W(\cdot+2^{-j}) - W]1_{[0,1-2^{-j}]}$:
\[
  \E\sup_{j\geq 1} \|\xi_j\|
  \lesssim\sup_{j\geq 1}\E\|\xi_j\|
  +\|(\sigma_j)_{j\geq 1}\|_{\Theta}.
\]
The first term is clearly smaller than $(\E\|W(1)\|^p)^{1/p}$.
By Lemma~\ref{lem:varianceW} and Example~\ref{ex:geometricSeq}, the
Orlicz norm can be computed as
\[\begin{aligned}
  \|(\sigma_j)_{j\geq 1}\|_{\Theta}
  \eqsim\|i_W\|\,\|(2^{-j/p})_{j\geq 1}\|_{\Theta}
  & \eqsim\|i_W\|\sqrt{\log(1-2^{-1/p})^{-1}}
  \\ & \eqsim (1+\sqrt{\log p})\|i_W\|.
\end{aligned}\] By Lemma~\ref{lem:extra}, this is smaller than
$(\E\|W(1)\|^p)^{1/p}$; indeed, it is much smaller when
$p\to\infty$. Thus, just like in Example~\ref{ex:mVsSigma}, we are
in a situation where the $m$ term totally dominates in the
estimate~\eqref{eq:estlinfty}. The proof of \eqref{eq:expBesov} is
complete.




Next, we show \eqref{eq:expBesovOrlicz}. The lower estimate follows
trivially from \eqref{eq:expBesov}.
For the upper estimate we write
\[\begin{split}
 \E &\|W\|_{B^{1/2}_{\Phi_2,\infty}(0,1;X)} \\
  &\leq \E\|W\|_{\mathfrak{L}^{\Phi_2}(0,1;X)}
  + \E\sup_{j\geq 1} 2^{j/2}\|W(\cdot+2^{-j})-W\|_{\mathfrak{L}^{\Phi_2}(0,1-2^{-j};X)}.
\end{split}\]

The first term can again be estimated using Doob's maximal
inequality, since
\[
  \E\|W\|_{\mathfrak{L}^{\Phi_2}(0,1;X)}\leq \E \|W\|_{L^\infty(0,1;X)}.
\]

The second term can be treated using Proposition \ref{prop:supestX}
with the $\mathfrak{L}^{\Phi_2}(0,1;X)$-valued Gaussian random
variables $\xi_j=2^{j/2}[W(\cdot+2^{-j}) - W]1_{[0,1-2^{-j}]}$.
Combining that Proposition with Remark~\ref{rem:expVsLp}, we have
\[
  \E\sup_{j\geq 1} \|\xi_j\|
  \lesssim\sup_{j\geq 1}\E\|\xi_j\|+\Big(\sum_{j\geq 1}\sigma_j^4\Big)^{1/4}.
\]
From Corollary~\ref{cor:varianceW} we get
\[
  \sigma_j\lesssim (\log 2^j)^{-1/2}\|i_W\|\eqsim j^{-1/2}\|i_W\|,
\]
so that the series sums up to $\lesssim\|i_W\|\lesssim\E\|W(1)\|$.

We then estimate $\E\|\xi_j\|$. 
By \eqref{eq:orliczineqL}, we have
\[\begin{aligned}
  \|f\|_{\mathfrak{L}^{\Phi_2}(0,1-2^{-j};X)} & \leq \|f\|_{\mathfrak{L}^{\Phi_2}(0,1;X)}
  \\ & \leq \|f\|_{L^{\Phi_2}(0,1;X)}
  = \inf_{\lambda>0} \frac{1}{\lambda}\int\limits_0^1 \exp(\lambda^2 \|f(t)\|^2) \, dt.
\end{aligned}\]

Therefore,
\[\begin{split}
  \E\|\xi_j\|
  &\leq \inf_{\lambda>0} \frac{1}{\lambda}
   \int\limits_0^1 \E \exp(\lambda^2 2^{j}\|W(t+2^{-j}) - W(t)\|^2) \, dt\\
  &= \inf_{\lambda>0} \frac{1}{\lambda}  \E \exp(\lambda^2 \|W(1)\|^2).
\end{split}\]
This may be estimated by expanding into a power series:
\[\begin{split}
  \frac{1}{\lambda} &\sum_{k\geq 0}\frac{\lambda^{2k}}{k!}\E\|W(1)\|^{2k} \\
  &\leq \frac{1}{\lambda}\Big[1+\sum_{k\geq 1}\frac{\lambda^{2k}}{k!}(K\sqrt{2k}\,\E\|W(1)\|)^{2k}\Big] \\
  &\leq\frac{1}{\lambda}\Big[1+\sum_{k\geq 1}(2e[\lambda K\E\|W(1)\|]^2)^k\Big],
\end{split}\]
where $K$ is an absolute constant from the Gaussian norm comparison
result \cite[Corollary~3.2]{LeTa}, and we used $k^k/k!\leq e^k$.
With the choice $\lambda=(2eK\E\|W(1)\|)^{-1}$, we find that
$\E\|\xi_j\|\lesssim\E\|W(1)\|$.
\end{proof}

{\em Acknowledgement} -- The authors thank Jan van Neerven for some
helpful comments

\def\cprime{$'$}
\providecommand{\bysame}{\leavevmode\hbox
to3em{\hrulefill}\thinspace}

\end{document}